\documentclass[12pt]{article}
\usepackage{amssymb}
\usepackage{epsfig}
\begin{document}
\renewcommand{\appendix}[1]
{\addtocounter{section}{1} \setcounter{equation}{0}
  \renewcommand{\thesection}{\Alph{section}}
  \section*{Appendix \thesection\protect\indent #1}
  \addcontentsline{toc}{section}{Appendix \thesection\ \ \ #1}
  }
\newcommand{\Section}[1]{\setcounter{equation}{0}\section{#1}}
\renewcommand{\theequation}{\thesection.\arabic{equation}}
\def\beq{\begin{equation}}
\def\eeq{\end{equation}}
\def\bea{\begin{eqnarray}}
\def\eea{\end{eqnarray}}
\def\A{{\cal A}}
\def\eps{\varepsilon}
\def\ha{\frac{1}{2}}
\def\hR{\hat R}
\def\bx{\bar x}
\def\by{\bar y}
\def\l{\lambda}
\def\hh{{1 \over 2}}
\def\iq{\frac{1}{q}}
\def\iqt{\frac{1}{q^2}}
\def\no{\nonumber\\}
\newcommand{\zet}{\hbox{$\mathbb Z$}}
\newcommand{\reals}{\hbox{$\mathbb R$}}
\newcommand{\la}{\langle}
\newcommand{\ra}{\rangle}
\newcommand{\poi}{Poincar\'e algebra}
\newcommand{\qq}{q^2+1}
\newcommand{\ci}{\! \circ \!}
\newcommand{\pn}{\par\noindent}
\thispagestyle{empty}
\hspace*{\fill}LMU-TPW 98--02\\
\hspace*{\fill}MPI-PhT/98--09\\
\hspace*{\fill}January 1998\\[3ex]
\begin{center}
{\Large\bf \boldmath$q$-Deformed Minkowski Space based on}\\
\end{center}
\begin{center}
{\Large\bf a \boldmath$q$-Lorentz Algebra}\\
\end{center}
\vspace{1cm}
\begin{center}
{\large B. L. Cerchiai and J. Wess}\\
\bigskip
{\it Sektion Physik der Ludwig-Maximilians-Universit\"at\\
Theresienstr. 37, D-80333 M\"unchen\\ and\\
Max-Planck-Institut f\"ur Physik\\
(Werner-Heisenberg-Institut)\\
F\"ohringer Ring 6, D-80805 M\"unchen\\
}
\end{center}
\vskip2cm
\centerline{\bf Abstract}
\bigskip
The Hilbert space representations of a non-commutative
$q$-deformed Minkowski space, its momenta and its Lorentz boosts
are constructed.

The spectrum of the diagonalizable space elements shows a lattice-like
structure with accumulation points on the light-cone.
\setcounter{page}{0}
\newpage
\section*{Introduction}

A non-commutative space-time structure emerges from quantum group
considerations.

More precisely, if we demand that space-time variables are modules
or co-modules of the $q$-deformed Lorentz group, then  they satisfy
commutation relations that make them elements of a non-comutative
space. The action of momenta on this space is non-commuative as well. 
The full structure is determined by the (co-)module property.

This algebra has been constructed in ref. \cite{lorek}.
It can serve as an explicit example of a non-commutative structure
for space-time.

This has the advantages that the $q$-deformed Lorentz group plays the
role of a cinematical group and thus determines many of the properties of
this space and allows explicit calculations. We have explicitly constructed
Hilbert space representations of the algebra and find that the
vectors in the Hilbert space can be determined by measuring the
time, the three-dimensional distance, the $q$-deformed angular momentum
and its third component. The eigenvalues of these observables form a
$q$-lattice with accumulation points on the light-cone. In a way physics
on the light-cone is best approximated by this $q$-deformation.
It is an interesting result that time-like and space-like regions serve
as basis for irreducible representations independently.
It will be shown however in a forthcoming paper \cite{forthcoming}
that these representations are linked together if we demand that the
observables are essentially selfadjoint operators.

The paper is organized as follows. We first present the algebra.
In chapter 2 we give explicit formulas for the matrix elements of
the elements of the algebra. This is the main result of this work and
can serve as a starting point for further investigations. In the
following chapters we give a rather detailed guide how these results
can be obtained, first for the space-time algebra (chapter 3), then
for the Lorentz algebra (chapter 4).

The algebra represented that far is isomorphic to the $q$-deformed
Poincar{\'e} algebra. We would just have to replace $X$ by $P$ to obtain
the respective representations \cite{pillin}.

In the next chapter (chapter 5) we enlarge the algebra by a scaling operator
and we introduce a canonical notation for labeling the states.

Finally in chapter 6 we construct the representations of the momenta
in the $X$-basis. There we learn that the full algebra cannot be represented
on the light-cone by itself. The points on the light-cone are limiting
points from the time-like and space-like regions.

\newpage

\Section{The algebra}

The algebra derived in ref. \cite{lorek} is generated by the elements $X^a$, 
(coordinates), $P^a$ (momenta), $V^{ab}$ ($q$-Lorentz transformations),
and $\Lambda^{\frac{1}{2}}$  (scaling operator).
\footnote{Capital letters ${}^A$ denote the three space indices $(+,-,3)$,
small letters ${}^a$ denote the four Minkowski indices $(+,-,3,0)$.
$\eps_{CBA}$ is the $q$-deformed $\eps$-tensor and $g_{AB}$ the Euclidean
metric, $\eta_{ab}$ the Lorentz metric. For the scalar product we write 
$X \ci Y= g_{AB} X^A Y^B$ (see also appendix A).}

The space is non-commutative:
\bea
\eps_{CB}{}^A X^B X^C &=&(1-q^2) X^0 X^A \label{XX} \\
X^0 X^C &=&X^C X^0
\nonumber
\eea
The momenta are subject to the same relations:
\bea
\eps_{CB}{}^A P^B P^C &=&(1-q^2) P^0 P^A \\
P^0 P^C &=&P^C P^0
\nonumber
\eea
The defining relations of the $q$-Lorentz algebra, as it acts on coordinates
and momenta, are more easily expressed in the ``Pauli''-notation:
\bea
V^{A0} &=& R^A + q^2 S^A \no
V^{0A} &=& -q^2 R^A -  S^A \\
V^{AB} &=& \eps_C{}^{AB} (R^C - S^C) \no
V^{00} &=& 0 \nonumber
\eea
In addition, we introduce an element $U$ that is related to
the Casimir operators of the $q$-Lorentz algebra:
\beq
U^2=1+{1 \over 2}(q^4-1)^2 (R \ci R + S \ci S)
\eeq
The $q$-Lorentz algebra:
\footnote{The $\hat{R}$ matrices are also defined in appendix A.}
\bea
\varepsilon_{CB}{}^A R^B R^C &=&\frac{1}{1 + q^2} UR^A \no
\varepsilon_{CB}{}^A S^B S^C &=& -\frac{1}{1 + q^2} US^A \label{RS} \\
R^A S^B &=& q^2 \hat{R}^{AB}{}_{CD} S^C R^D \nonumber
\eea
\[
UR^A=R^A U , \qquad  U S^A =S^A U  
\]
The coordinates ``transform'' under the $q$-Lorentz transformations:
\bea
R^AX^0 &=& \frac{1}{q}  \frac{q^4 + 1}{q^2 + 1} X^0 R^A  + \frac{1}{q}
\frac{q^2 - 1}{q^2 + 1}
 \varepsilon_{LM}{}^A X^M R^L- \frac{q}{(1 + q^2)^2} X^A U\nonumber\\
R^AX^B &=& \frac{1}{1 + q^2}  \Big[q( 1 + q^2) X^A R^B  - \frac{1}{q} (q^2 - 1)
 \varepsilon_C{}^ {AB} X^0 R^C \no
&-& \frac{1}{q} (q^2 - 1) g^{AB} g_{MC} X^M R^C - \frac{2}{q}
\varepsilon^{ABG} \varepsilon_{STG} X^T R^S\nonumber\\
&-& \frac{1}{q} \frac{1}{1+q^2} g^{AB} X^0 U +  \frac{1}{q} \frac{1}{1+q^2}
\varepsilon_M{}^{AB} X^M U \Big]\no[5mm]
S^AX^0 &=& \frac{1}{q}  \frac{q^4 + 1}{q^2 + 1} X^0 S^A  + \frac{1}{q}
\frac{q^2 - 1}{q^2 + 1}
 \varepsilon_{LM}{}^A X^M S^L- \frac{1}{q (1 + q^2)^2} X^A U\nonumber\\
S^AX^B &=& \frac{1}{1 + q^2}  \Big[ \frac{1}{q} ( 1 + q^2) X^A S^B 
- \frac{1}{q} (q^2 - 1)
 \varepsilon_C{}^ {AB} X^0 S^C \label{SX} \\
&+& q (q^2 - 1) g^{AB} g_{MC} X^M S^C - \frac{2}{q}
\varepsilon^{ABG} \varepsilon_{STG} X^T S^S\nonumber\\
&-& \frac{q}{1+q^2} g^{AB} X^0 U - \frac{1}{q} \frac{1}{1+q^2}
\varepsilon_M{}^{AB}  X^M U \Big] \no[5mm]
U X^0&=&{1 \over q} {{q^4+1} \over {q^2+1}} X^0 U-{1 \over q} (q^2-1)^2
X \circ R \no
U X^A&=&{1 \over q} {{q^4+1} \over {q^2+1}} X^A U -q(q^2-1)^2 X^0 R^A
-\iq (q^2-1)^2 \varepsilon_{CB}{}^A X^B R^C \nonumber
\eea

The momenta have the same transformation law. 

The scaling operator acts as follows:
\bea
 \Lambda^{{1 \over 2}} X^a &=& \iq X^a \Lambda^{{1 \over 2}} \no
 \Lambda^{{1 \over 2}} P^a &=& q P^a  \Lambda^{{1 \over 2}}
\label{la} \\
 \Lambda^{{1 \over 2}} V^{ab} &=& V^{ab} \Lambda^{{1 \over 2}} \no
\Lambda^{{1 \over 2}} U  &=& U \Lambda^{{1 \over 2}} \nonumber
\eea 

The relations that generalizes the Heisenberg commutation relations are:
${}^2$
\beq
\begin{array}{c}
P^a X^b - q^{-2} \hat{R}_{II}^{-1}{}^{ab}{}_{cd} X^c P^d = \\[1mm]
- \frac{i}{2} \Lambda^{-{1 \over 2}} \left\{(1 + q^4) \eta^{ab} U 
+ q^2 (1 - q^4) V^{ab} \right\}
\end{array}
\label{PX}
\eeq
The $q$-Heisenberg algebra (\ref{PX}) does not separate from the $q$-Lorentz
algebra for $q \neq 1$. The relation (\ref{PX}) tells us how
to commute $X^a$ and $P^b$ and how to define orbital angular momentum
in terms of the space and momentum operators. It is not possible to
define $V^{ab}$ in terms of an $X,P$ ordered expression.

{From} orbital angular momentum we expect additional relations -- the orbital
angular momentum is orthogonal to the coordinates and momenta. These 
relations follows from the defining relations of our algebra and they are:
\bea
g_{AB}  X^A (R^B-q^2 S^B)&=0
\label{extra3} \\
X^0(S^A-q^2 R^A)-\varepsilon_{CB}{}^A X^B(R^C+S^C)&=0 \nonumber
\eea
The same with $X$  replaced by $P$.
Not all representations of the Lorentz group can be realized as angular
momentum. We expect a relation for the Casimir operators. It follows from
the algebra that:
\beq
R \ci R = S \ci S
\eeq
For the physical interpretation and for the representations of this algebra
the conjugation properties are very important. They are:
\bea
\overline{X^0} = X^0  &,& \quad   \overline{X^A} = g_{AB} X^B
\no
\overline{P^0} = P^0 & , &\quad \overline{P^A} = g_{AB} P^B
\no
\overline{R^A} = -g_{AB} S^B & , &\quad \overline{S^A} = -g_{AB} R^B 
\label{star} \\
\overline{U} &=& U
\no
\overline{\Lambda^{1/2}} &=& q^4 \Lambda^{-1/2}
\nonumber
\eea
These conjugation relations are consistent with the algebra.

Finally, we identify the three-dimensional rotations in the algebra.
They have to commute with $X^0$ and $X \ci X$ Such operators have been
found in \cite{cerchiai} and they are:
\beq
L^A = {{q^2+1} \over q^2}(U S^A-U R^A+(q^4-1) \eps_{CB}{}^A R^B S^C)
\eeq
They commute with $X^0$ and $P^0$ as well as with all ``scalars'' in our
algebra formed with the metric $g_{AB}$, such as $X \ci X$, $P \ci P$,
$R \ci R$, $S \ci R$, etc.

To write the $L$ algebra in a familiar way we define an additional
element:
\beq
W=U^2-q^2(q^4-1)^2 R \ci S
\eeq
and find:
\bea
\varepsilon_{BC}{}^A L^C L^B&=&-{W \over q^2} L^A \\
q^4(q^2-1)^2L\circ L&=&W^2-1
\nonumber
\eea

The $SU_q(2)$ algebra was written in this form in \cite{symm,witten}. We 
identify the $SU_q(2)$ generators:
\bea
T^+ &=& q^2 \sqrt{1 + q^2}  \quad \tau^{1/2} L^+ 
\no
T^- &=&-q^3 \sqrt{1 + q^2}  \quad \tau^{1/2} L^- \label{mapLT} \\
\tau_3 &=& \tau \nonumber
\eea
with $\tau^{-\frac{1}{2}}=W + q^2 (1 - q^2) L^3$.

The $T$ algebra is the familiar one:
\bea
\iq T^+ T^- -q T^- T^+&=&{1-\tau_3 \over q-\iq} \no
\tau_3 T^+&=&{1 \over q^4} T^+ \tau_3 \label{TT}\\
\tau_3 T^-&=&q^4 T^- \tau_3 \nonumber
\eea
Its Casimir operator is:
\beq
\vec{T}^2=q T^- T^+ +{q \over (q-\iq)^2} \tau_3^{-\frac{1}{2}} + \iq
{1 \over (q-\iq)^2} \left(\tau_3^{\frac{1}{2}} -q^2-1 \right)
\label{casimir}
\eeq
The conjugation properties are:
\bea
\overline{T^+}=\iqt T^-, & \quad \overline{T^-}=q^2 T^+ 
\label{Tconj} \\
\overline{\tau_3}=\tau_3 \nonumber
\eea
The vectors $X^A$ and $P^A$ transform as follows:
\bea
L^A X^B &=& g^{AB} X \circ L 
-\frac{1}{q^2}\eps_{KC}{}^{A} \varepsilon_{D}{}^{KB} X^C L^D 
-\frac{1}{q^4}\varepsilon_{C}{}^{AB}X^C W
\label{LX} \\
WX^A &=&(q^2+\frac{1}{q^2}-1)X^A W +(q^2-1)^2 \varepsilon_{DC}{}^{A} X^C L^D
\no
W X^0 &=& X^0 W
\nonumber
\eea
Now we have all the relations that allow us to study the representations of 
this algebra.

A complete set of commuting operators is $X^0$, $X \ci X$, $\vec{T}^2$
and $\tau_3$.

\Section{The matrix elements}

In this chapter we present the matrix elements of all the members of the
algebra.

The states are labeled by the quantum numbers $j$, $m$, $n$ and $M$.
The quantum numbers $j$ and $m$ refer to the $q$-deformed angular
momentum. The quantum numbers $n$ and $M$ label the eigenvalues of
the time $X^0$ and the three-dimensional radius $X \ci X$.

There are inequivalent representations for the time-like and
space-like regions.

\pn Space-like: $s^2= t^2-r^2<0$:
\begin{eqnarray*} 
M=-\infty \ldots \infty \\
n=-\infty \ldots \infty \\
j=0 \ldots \infty
\end{eqnarray*} 
\bea
X^0 \; |j,m,n,M \ra &=& {l_0 q^M \over [2]} \lambda [n] \; |j,m,n,M \ra \\
X \ci X \; |j,m,n,M \ra &=&  {l_0^2 q^{2M} \over [2]^2} \{n+1\} \{n-1\} \; 
|j,m,n,M \ra \nonumber
%\vec{T}^2 \; |j,m,n,M \ra &=& [j] [j+1] \; |j,m,n,M \ra \no
%\tau_3   \; |j,m,n,M \ra &=& q^{-4m} \; |j,m,n,M \ra \nonumber
\eea
\pn Time-like: $s^2= t^2-r^2>0$:
\begin{eqnarray*} 
M=-\infty \ldots \infty \\
n=0 \ldots \infty \\
j=0 \ldots n
\end{eqnarray*} 
\bea
X^0 \; |j,m,n,M \ra &=& {t_0 q^M \over [2]} \{n+1\} \; |j,m,n,M \ra \\
X \ci X \; |j,m,n,M \ra &=&   {t_0^2 q^{2M} \lambda^2 \over [2]^2} [n+2] [n] 
\; |j,m,n,M \ra \nonumber
%\vec{T}^2 \; |j,m,n,M \ra &=& [j] [j+1] \; |j,m,n,M \ra \no
%\tau_3   \; |j,m,n,M \ra &=& q^{-4m} \; |j,m,n,M \ra \nonumber
\eea

We use the notation throughout this paper:
\bea
[a]&=&{q^a -q^{-a} \over q-q^{-1}} \\
\{a\}&=&q^a +q^{-a}
\eea
and
\beq
\lambda=q-\iq
\eeq
The spectrum of the operators $X^0$ and $X \ci X$ is shown in fig. 
The parameters $|t_0|$, $l_0$ range from $1$ to $q$ and label inequivalent 
representations. $t_0$ can be positive (forward cone) and negative 
(backward cone).

The states are orthonormal: 
\beq
\la j',m',n',M' | j,m,n,M \ra=\delta_{j',j} \delta_{m',m}\delta_{n',n}
\delta_{M',M}
\eeq
The matrix elements of $X^A$, $R^A$, $S^A$ and $P^A$ can
be expressed in terms of reduced matrix elements. The explicit formulas
are given in (\ref{redXm}).

\medskip

\pn {\bf Space-like:}

\medskip

\pn{\it Reduced matrix elements of $X^-$:}

\[
\begin{array}{lcll}
\la j+1,n,M \| X^- \| j,n,M \ra&=& l_0 q^{M+j}
{\sqrt{\displaystyle \{n+j+1\}\{n-j-1\}}
\over \displaystyle \{j+1\} \sqrt{[2][2j+1] [2j+3]}} 
& \quad \mbox{ for } j \ge 0\\[8mm]
\la j,n,M \| X^- \| j,n,M \ra &=&
-q^{-1}  {\displaystyle l_0 q^M [n]\lambda^2 \over 
\displaystyle \sqrt{[2]} \{j\}\{j+1\}} 
& \quad \mbox{ for } j \ge 1 \\[8mm]
\la j,n,M \| X^- \| j+1,n,M \ra &=& -l_0 q^{M-j-2}
{\sqrt{\displaystyle \{n+j+1\}\{n-j-1\}}
\over \displaystyle \{j+1\} \sqrt{[2][2j+1] [2j+3]}}
& \quad \mbox{ for } j \ge 0
\end{array}
\]

\pn {\it Reduced matrix elements of $R^-$:}

\[
\begin{array}{lcll}
\lefteqn{\la j+1,n',M\|R^-\|j,n,M \ra 
=(\delta_{n',n+1}+\delta_{n',n-1})
{(n'-n)q^{2j-1} \over  \{j+1\} [2]^{\frac{3}{2}} \lambda \sqrt{[2j+1][2j+3]}}}
\\[3mm]
& \cdot & {\displaystyle \sqrt{\{(n'-n)(j+1)+n'\}\{(n'-n)j+n'\} 
\over \{n\}\{n'\}}} 
& \quad \mbox{ for } j \ge 0 \makebox{\hspace*{3.8cm}} \\[8mm]
\end{array}
\]
\[
\begin{array}{lcll}
\lefteqn{\la j,n',M\|R^-\|j,n,M \ra=(\delta_{n',n+1}+\delta_{n',n-1})
{\displaystyle q^{-3} \over \displaystyle \{j+1\} \{j\} [2]^{\frac{3}{2}}}}
\\[3mm]
&\cdot &{\displaystyle \sqrt{\{(n'-n)j+n'\}\{n-(n'-n)j\} \over \{n \}\{n'\}}}
& \quad \mbox{ for $j \ge 1$} \makebox{\hspace*{3.8cm}} \\[8mm]
\lefteqn{\la j,n',M\|R^-\|j+1,n,M \ra =-(\delta_{n',n+1}+\delta_{n',n-1})
{(n'-n)q^{-2j-5} \over  \{j+1\} [2]^{\frac{3}{2}} \lambda \sqrt{[2j+1][2j+3]}}}
 \\[3mm]
& \cdot & {\displaystyle \sqrt{\{n'-(n'-n)(j+1)\}\{n-(n'-n)(j+1)\} 
\over \{n\}\{n'\}}}
& \quad \mbox{ for $j \ge 0$} \makebox{\hspace*{3.8cm}}
\end{array}
\]

\newpage

\pn {\it Reduced matrix elements of $S^-$:}

\[
\begin{array}{lcll}
\lefteqn{\la j+1,n',M\|S^-\|j,n,M \ra =(\delta_{n',n+1}+\delta_{n',n-1})
{(n'-n)q^{-3} \over  \{j+1\} [2]^{\frac{3}{2}} \lambda
\sqrt{[2j+1][2j+3]}}} \\[3mm]
&\cdot &{\displaystyle \sqrt{\{(n'-n)(j+1)+n'\}\{(n'-n)j+n'\} 
\over \{n\}\{n'\}}} 
& \quad \makebox{ for $j \ge 0$} \makebox{\hspace*{4cm}} \\[8mm]
\lefteqn{\la j,n',M\|S^-\|j,n,M \ra =-(\delta_{n',n+1}+\delta_{n',n-1})
{\displaystyle q^{-3} \over \displaystyle \{j+1\} \{j\} [2]^{\frac{3}{2}}}}
\\[3mm]
& \cdot & {\displaystyle \sqrt{\{(n'-n)j+n'\}\{n-(n'-n)j\} \over \{n \}\{n'\}}}
& \quad \makebox{ for $j \ge 1$} \makebox{\hspace*{4cm}} \\[8mm]
\lefteqn{\la j,n',M\|S^-\|j+1,n,M \ra =-(\delta_{n',n+1}+\delta_{n',n-1})
{(n'-n)q^{-3} \over  \{j+1\} [2]^{\frac{3}{2}} \lambda 
\sqrt{[2j+1][2j+3]}}} \\[3mm]
& \cdot & {\displaystyle \sqrt{\{n'-(n'-n)(j+1)\}\{n-(n'-n)(j+1)\} 
\over \{n\}\{n'\}}}
& \quad \makebox{ for $j \ge 0$} \makebox{\hspace*{4cm}}
\end{array}
\]

\pn {\it Reduced matrix elements of $P^-$:}

\[
\begin{array}{lcll}
\lefteqn{\la j+1,n+1,M' \| P^- \| j,n,M \ra={\displaystyle i \over 
\displaystyle  2}
(\delta_{M',M+1} q^{2+2j-n}-\delta_{M',M-1} q^{2+n})} \\[3mm]
& \cdot &
{\displaystyle 1 \over \displaystyle \{j+1\} \lambda l_0 q^M}
 \sqrt{{ \displaystyle[2] \over \displaystyle [2j+1][2j+3]}}
\sqrt{\displaystyle \{n+j+2\}\{n+j+1\} \over \displaystyle \{n\}\{n+1\}}
& \quad \makebox{ for $j \ge 0$} \makebox{\hspace*{1.5cm}}
\\[8mm]
\lefteqn{\la j+1,n-1,M' \| P^- \| j,n,M \ra={\displaystyle i \over 
\displaystyle  2} (\delta_{M',M+1} q^{2+2j+n}-\delta_{M',M-1} q^{2-n})} \\[3mm]
& \cdot & {\displaystyle 1 \over \displaystyle \{j+1\} \lambda l_0 q^M}
\sqrt{{\displaystyle [2] \over\displaystyle  [2j+1][2j+3]}}
\sqrt{{\displaystyle \{n-j-2\}\{n-j-1\} \over \displaystyle \{n\}\{n-1\}}}
& \quad \makebox{ for $j \ge 0$} \makebox{\hspace*{1.5cm}} \\[8mm]
\lefteqn{\la j,n+1,M' \| P^- \| j,n,M \ra={\displaystyle i \lambda \over 
\displaystyle  2\{j\} \{j+1\} \sqrt{[2]} l_0 q^M} 
(\delta_{M',M+1} q^{-n}+\delta_{M',M-1} q^{2+n})} \\[3mm]
& \cdot & \sqrt{\displaystyle \{n+j+1\}\{n-j\} \over \displaystyle 
\{n\}\{n+1\}} & \quad \mbox{ for } j \ge 1 \makebox{\hspace*{1.5cm}} \\[8mm]
\end{array}
\]
\[
\begin{array}{lcll}
\lefteqn{\la j,n-1,M' \| P^- \| j,n,M \ra=-{\displaystyle i \lambda \over 
\displaystyle  2 \{j\} \{j+1\} \sqrt{[2]}l_0 q^M} 
(\delta_{M',M+1} q^n+\delta_{M',M-1} q^{2-n})} \\[3mm]
& \cdot & \sqrt{{\displaystyle \{n-j-1\}\{n+j\} \over \displaystyle 
\{n\}\{n-1\}}} & \quad \mbox{ for } j \ge 0 \makebox{\hspace*{1.5cm}} \\[8mm]
\lefteqn{\la j,n+1,M' \| P^- \| j+1,n,M \ra=-{\displaystyle i \over 
\displaystyle  2}
(\delta_{M',M+1} q^{-2-2j-n}-\delta_{M',M-1} q^{2+n})} \\[3mm]
& \cdot & {\displaystyle 1 \over \displaystyle \{j+1\} \lambda l_0 q^M} 
 \sqrt{{ \displaystyle[2] \over \displaystyle [2j+1][2j+3]}}
\sqrt{\displaystyle \{n-j\}\{n-j-1\} \over \displaystyle \{n\}\{n+1\}} 
& \quad \mbox{ for } j \ge 0 \makebox{\hspace*{1.5cm}} \\[8mm]
\lefteqn{\la j,n-1,M' \| P^- \| j+1,n,M \ra=-{\displaystyle i \over 
\displaystyle  2} 
(\delta_{M',M+1} q^{-2-2j+n}-\delta_{M',M-1} q^{2-n})} \\[3mm]
& \cdot & {\displaystyle 1 \over \displaystyle \{j+1\} \lambda l_0 q^M} 
 \sqrt{{\displaystyle [2] \over\displaystyle  [2j+1][2j+3]}}
\sqrt{{\displaystyle \{n+j+1\}\{n+j\} \over \displaystyle \{n\}\{n-1\}}}  
& \quad \mbox{ for } j \ge 0 \makebox{\hspace*{1.5cm}}
\end{array}
\]

\pn {\it Matrix elements of $P^0$:}

\[
\begin{array}{lcll}
\lefteqn{\la j,m,n+1,M' | P^0 | j,m,n,M \ra=-{i \over 2 \lambda l_0  q^M}
(\delta_{M',M+1} q^{1-n}+\delta_{M',M-1} q^{3+n})} \\[3mm]
& \cdot &{\displaystyle \sqrt{\{n-j\}\{n+j+1\} \over \{n\}\{n+1\}}} 
& \quad \mbox{for } j \ge 0 \makebox{\hspace*{8cm}} \\[8mm]
\lefteqn{\la j,m,n-1,M' | P^0 | j,m,n,M \ra={i \over 2 \lambda l_0  q^M}
(\delta_{M',M+1} q^{1+n}+\delta_{M',M-1} q^{3-n})} \\[3mm]
& \cdot & {\displaystyle \sqrt{\{n-j-1\}\{n+j\} \over \{n\}\{n-1\}}}
& \quad \mbox{for } j \ge 0 \makebox{\hspace*{8cm}}
\end{array}
\]

\pn {\it Matrix elements of $U$:}

\[
\begin{array}{lcll}
&\la j,m,n,M | U |j,m,n+1,M \ra=\la j,m,n+1,M| U |j,m,n,M\ra=
\makebox{\hspace*{1cm}} \\[3mm]
&{\displaystyle {1 \over [2]} \sqrt{\{n-j\}\{n+j+1\} \over 
\{n\}\{n+1\}}} 
& \quad \mbox{for } j \ge 0 
\end{array}
\]

\pn  {\it Matrix elements of $\Lambda$:}

\[
\la j,m,n,M+1 |\Lambda^{\frac{1}{2}} |j,m,n,M \ra = 
q^2 
\]

%\begin{eqnarray*}
%\la j, m+1, n',M' | A^+ | j, m, n,M  \ra &=&-q^{m+2} \sqrt{[j+m+1][j-m]} 
%\la j, n',M' \| A^- \| j, n,M \ra \\
%\la j+1, m+1, n',M' | A^+ \| j, m, n, M \ra &=&q^{m-2j} 
%\sqrt{[j+m+1][j+m+2]} \\
%&&\cdot \la j+1, n',M' \| A^- \| j, n,M \ra  \\
%\la j-1, m+1, n',M' | A^+ | j, m, n,M \ra &=&q^{m+2j+2}
%\sqrt{[j-m][j-m-1]} \\
%&& \cdot \la j-1, n',M' \| A^- \| j, n,M  \ra
%\end{eqnarray*}
%\bea
%\la j, m-1, n', M' | A^- | j, m, n, M\ra &=&q^m
%\sqrt{[j+m][j-m+1]} \la j, n', M' \| A^- \| j, n, M \ra  \no
%\la j+1, m-1, n', M' | A^- | j, m, n, M\ra &=&q^m
%\sqrt{[j-m+1][j-m+2]} \label{redAm}\\
%&&\cdot \la j+1, n', M' \| A^- \| j, n, M \ra  \no
%\la j-1, m-1, n', M' | A^- | j, m, n, M\ra &=&q^m
%\sqrt{[j+m][j+m-1]} \la j-1, n', M' \| A^- \| j, n, M \ra  \nonumber
%\eea
%\begin{eqnarray*}
%\la j, m, n', M' | A^3 | j, m, n, M\ra &=&q^{3 \over 2} 
%{{\sqrt{1+q^2}}\over{q^2-1}} \{q^{2m}- {{q^{2j+1}+q^{-2j-1}}
%\over{q+q^{-1}}}\} \la j, n', M' \| A^- \| j, n, M \ra  \\
%\la j+1, m, n', M' | A^3 | j, m, n, M\ra &=&q^{m-j-{1 \over 2}}
%\sqrt{1+q^2} \sqrt{[j-m+1][j+m+1]} \\
%&& \cdot \la j+1, n', M' \| A^- \| j, n, M \ra \no
%\la j-1, m, n', M' | A^3 | j, m, n, M\ra &=&-q^{m+j+{1 \over 2}}
%\sqrt{1+q^2} \sqrt{[j-m][j+m]} \\
%&& \cdot \la j-1, n', M' \| A^- \| j, n, M \ra
%\end{eqnarray*}
%where $\la j', n', M' \| A^- \| j, n, M \ra$ are the reduced matrix 
%elements, which do not depend on $m$.

\newpage

\pn {\bf Time-like:}

\medskip

\pn{\it Reduced matrix elements of $X^-$:}

\[
\begin{array}{lcll}
\la j+1,n,M \| X^- \| j,n,M \ra&=& t_0 q^{M+j} \lambda
{\sqrt{\displaystyle [n-j][n+j+2]}
\over \displaystyle \{j+1\} \sqrt{[2][2j+1] [2j+3]}} 
& \quad \mbox{ for } j \ge 0 \\[8mm]
\la j,n,M \| X^- \| j,n,M \ra &=& 
-q^{-1} \lambda {\displaystyle t_0 q^M \{n+1\} \over \displaystyle \sqrt{[2]}
\{j\}\{j+1\}} \quad \mbox{ for } j \ge 1 \\[8mm]
\la j,n,M \| X^- \| j+1,n,M \ra&=&-t_0 q^{M-j-2}\lambda
{\sqrt{\displaystyle [n-j][n+j+2]}
\over \displaystyle \{j+1\} \sqrt{[2][2j+1] [2j+3]}}
& \quad \mbox{ for } j \ge 0
\end{array}
\]

\pn {\it Reduced matrix elements of $R^-$:}

\[
\begin{array}{lcll}
\lefteqn{\la j+1,n',M\|R^-\|j,n,M \ra =(\delta_{n',n+1}+\delta_{n',n-1})
{(n'-n) q^{2j-1} \over  \{j+1\} [2]^{\frac{3}{2}} \lambda
\sqrt{[2j+1][2j+3]}}} \\[3mm]
& \cdot & {\displaystyle \sqrt{[(n'-n)(j+1)+n'+1][(n'-n)j+n'+1] 
\over [n+1][n'+1]}} 
& \quad \makebox{ for $j \ge 0$} 
\makebox{\hspace*{3cm}} \\[8mm]
\lefteqn{\la j,n',M\|R^-\|j,n,M \ra=(\delta_{n',n+1}+\delta_{n',n-1})
{\displaystyle q^{-3} \over \displaystyle \{j+1\} \{j\} [2]^{\frac{3}{2}}}}
\\[3mm]
& \cdot & {\displaystyle \sqrt{[(n'-n)j+n'+1][n-(n'-n)j+1] \over 
[n+1][n'+1]}} 
& \quad \makebox{ for $j \ge 1$} 
\makebox{\hspace*{3cm}}\\[8mm]
\lefteqn{\la j,n',M \|R^-\|j+1,n,M \ra =-(\delta_{n',n+1}+\delta_{n',n-1})
{(n'-n) q^{-2j-5} \over  \{j+1\} [2]^{\frac{3}{2}} \lambda 
\sqrt{[2j+1][2j+3]}}} \\[3mm]
& \cdot & {\displaystyle \sqrt{[n'-(n'-n)(j+1)+1][n-(n'-n)(j+1)+1] 
\over [n+1][n'+1]}} 
& \quad \makebox{ for $j \ge 0$} 
\makebox{\hspace*{3cm}}
\end{array}
\]

\pn {\it Reduced matrix elements of $S^-$:}

\[
\begin{array}{lcll}
\lefteqn{\la j+1,n',M\|S^-\|j,n,M \ra =(\delta_{n',n+1}+\delta_{n',n-1})
{(n'-n)q^{-3} \over  \{j+1\} [2]^{\frac{3}{2}} \lambda \sqrt{[2j+1][2j+3]}}}
 \\[3mm]
& \cdot & {\displaystyle \sqrt{[(n'-n)(j+1)+n'+1][(n'-n)j+n'+1] 
\over [n+1][n'+1]}} 
& \quad  \makebox{ for $j \ge 0$} 
\makebox{\hspace*{3.85cm}}
\end{array}
\]
\[
\begin{array}{lcll}
\lefteqn{\la j,n',M\|S^-\|j,n,M \ra=-(\delta_{n',n+1}+\delta_{n',n-1})
{\displaystyle q^{-3} \over \displaystyle \{j+1\} \{j\} [2]^{\frac{3}{2}}}}
\\[3mm]
& \cdot & {\displaystyle \sqrt{[(n'-n)j+n'+1][n-(n'-n)j+1] \over 
[n+1][n'+1]}}
& \quad  \makebox{ for $j \ge 1$} 
\makebox{\hspace*{3.50cm}} \\[8mm]
\lefteqn{\la j,n',M\|S^-\|j+1,n,M \ra =-(\delta_{n',n+1}+\delta_{n',n-1})
{(n'-n)q^{-3} \over  \{j+1\} [2]^{\frac{3}{2}} \lambda 
\sqrt{[2j+1][2j+3]}}} \\[3mm]
& \cdot & {\displaystyle \sqrt{[n'-(n'-n)(j+1)+1][n-(n'-n)(j+1)+1] 
\over [n+1][n'+1]}} 
& \quad \makebox{ for $j \ge 0$} 
\makebox{\hspace*{3.50cm}}
\end{array}
\]

\pn {\it Reduced matrix elements of $P^-$:}

\[
\begin{array}{lcll}
\lefteqn{\la j+1,n+1,M' \| P^- \| j,n,M \ra={\displaystyle i 
\over \displaystyle 2} (\delta_{M',M+1} q^{1+2j-n} +\delta_{M',M-1} q^{3+n})}
\\[3mm]
& \cdot & {\displaystyle 1 \over \displaystyle t_0 q^M \lambda \{j+1\}}
\sqrt{\displaystyle [2] \over \displaystyle [2j+1][2j+3]}
\sqrt{\displaystyle [n+j+3][n+j+2] \over \displaystyle[n+2][n+1]} 
& \quad \mbox{ for } j \ge 0 \makebox{\hspace*{2.8cm}}\\[8mm]
\lefteqn{\la j+1,n-1,M' \| P^- \| j,n,M \ra=
-{\displaystyle i \over \displaystyle 2}
(\delta_{M',M+1} q^{3+2j+n} +\delta_{M',M-1} q^{1-n})} \\[3mm]
& \cdot &{\displaystyle 1 \over \displaystyle t_0 q^M \lambda \{j+1\}}
 \sqrt{\displaystyle [2] \over \displaystyle [2j+1][2j+3]}
\sqrt{\displaystyle [n-j-1][n-j] \over \displaystyle [n][n+1]} 
& \quad \mbox{ for } j \ge 0 \\[8mm]
\lefteqn{\la j,n+1,M' \| P^- \| j,n,M \ra={\displaystyle i \lambda  
\over \displaystyle 2 \{j\} \{j+1\} \sqrt{[2]} t_0 q^M} 
(\delta_{M',M+1} q^{-1-n} -\delta_{M',M-1} q^{3+n})} \\[3mm]
& \cdot & \sqrt{\displaystyle [n+j+2][n-j+1] \over \displaystyle[n+2][n+1]}
\quad \mbox{ for } j \ge 1 \\[8mm]
\lefteqn{\la j,n-1,M' \| P^- \| j,n,M \ra={\displaystyle i \lambda
\over \displaystyle 2 \{j\} \{j+1\} \sqrt{[2]}t_0 q^M}
(\delta_{M',M+1} q^{1+n} -\delta_{M',M-1} q^{1-n})} \\[3mm] 
& \cdot & \sqrt{\displaystyle [n+j+1][n-j] \over \displaystyle [n][n+1]}
& \quad \mbox{ for } j \ge 1\\[8mm]
\end{array}
\]
\[
\begin{array}{lcll}
\lefteqn{\la j,n+1,M' \| P^- \| j+1,n,M \ra=
-{\displaystyle i\over \displaystyle 2}
(\delta_{M',M+1} q^{-3-2j-n} +\delta_{M',M-1} q^{3+n})} \\[3mm]
& \cdot & {\displaystyle 1 \over \displaystyle \lambda \{j+1\}t_0 q^M}
 \sqrt{\displaystyle [2] \over \displaystyle [2j+1][2j+3]}
\sqrt{\displaystyle [n-j+1][n-j] \over \displaystyle[n+2][n+1]} 
& \quad \mbox{ for } j \ge 0 \\[8mm]
\lefteqn{\la j,n-1,M' \| P^- \| j+1,n,M \ra={\displaystyle i 
\over \displaystyle 2}
(\delta_{M',M+1} q^{-1-2j+n} +\delta_{M',M-1} q^{1-n})} \\[3mm]
& \cdot & {\displaystyle 1 \over \displaystyle \lambda \{j+1\}t_0 q^M }
 \sqrt{\displaystyle [2] \over \displaystyle [2j+1][2j+3]}
\sqrt{\displaystyle [n+j+2][n+j+1] \over \displaystyle [n][n+1]} 
& \quad \mbox{ for } j \ge 0
\end{array}
\]

\pn {\it Matrix elements of $P^0$:}

\[
\begin{array}{lcll}
\la j,m,n+1,M' | P^0 | j,m,n,M \ra&=&-{i \over 2 \lambda t_0 q^M}
(\delta_{M',M+1} q^{-n}-\delta_{M',M-1} q^{4+n}) \\[3mm]
\cdot {\displaystyle \sqrt{[n-j+1][n+j+2] \over [n+1][n+2]}} 
&&\quad \mbox{for } j \ge 0\\[8mm]
 \la j,m,n-1,M' | P^0 | j,m,n,M \ra&=&-{i \over 2 \lambda t_0 q^M}
(\delta_{M',M+1} q^{2+n}-\delta_{M',M-1} q^{2-n})\\[3mm]
\cdot {\displaystyle \sqrt{[n-j][n+j+1] \over [n][n+1]}}
\end{array}
\]

\pn {\it Matrix elements of $U$:}

\[
\begin{array}{lcll}
&\la j,m,n,M| U |j,m,n+1,M \ra=\la j,m,n+1,M| U |j,m,n,M\ra=\\[3mm]
&{\displaystyle{1 \over [2]} 
{\displaystyle \sqrt{[n-j+1][n+j+2] \over [n+1][n+2]}}}
\end{array}
\]

\pn  {\it Matrix elements of $\Lambda$:}

\beq
\la j,m,n,M+1 |\Lambda^{\frac{1}{2}} |j,m,n,M \ra = q^2 
\eeq

\centerline{\epsfig{figure=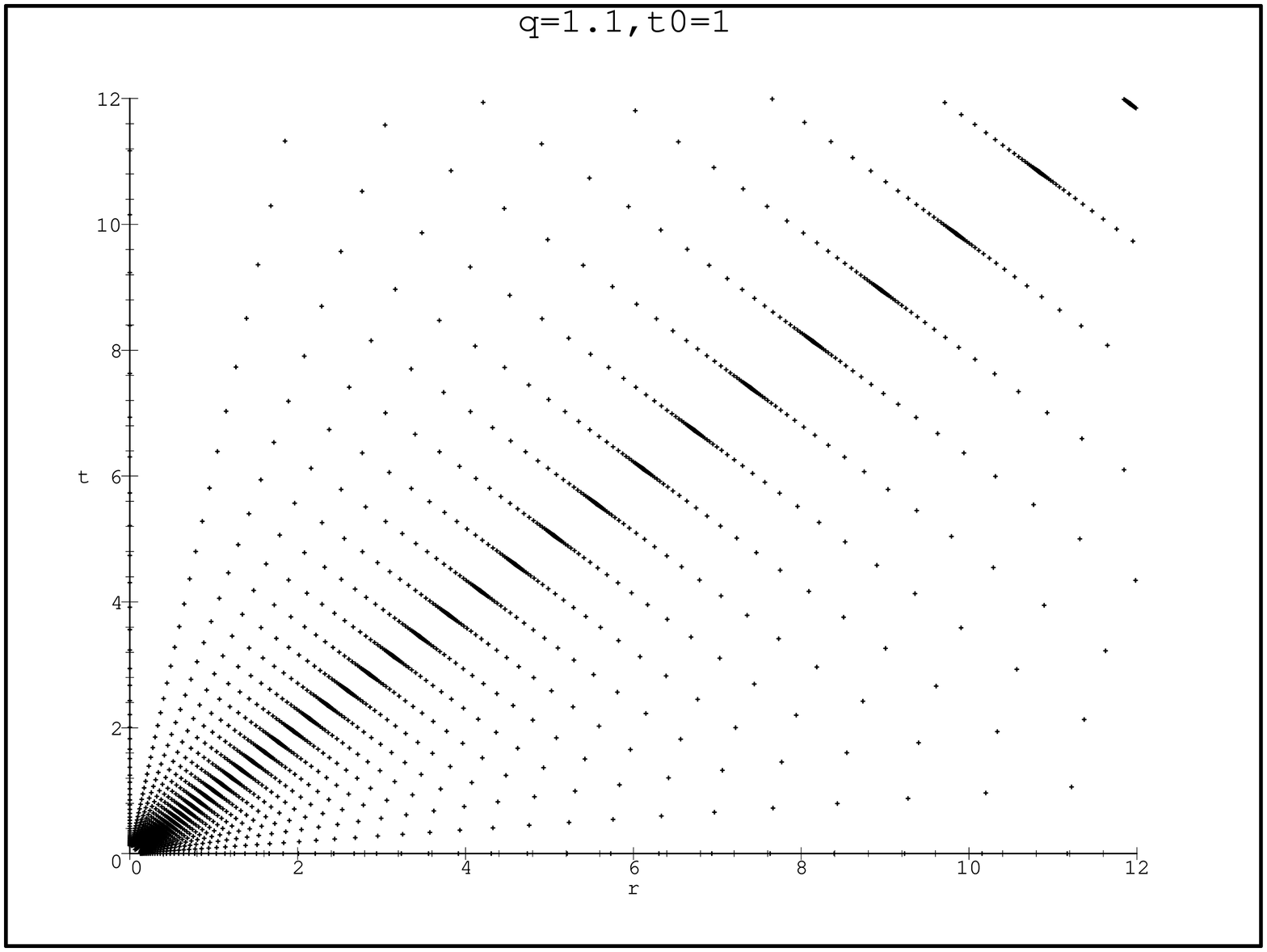,height=175mm,width=175mm,angle=270}}
\begin{figure}[ht]
\vspace{-1.5cm}
\caption{Admissible values of $t$ versus those of $r$ for $q=1.1$ and $t_0=1$.}
\end{figure}

\newpage
\Section{Matrix elements of the coordinates}

In this chapter we indicate how to construct the matrix elements of the 
coordinates $X$.

We assume $X^0$ and $X \ci X$, as well as the elements $\vec{T}^2$ and $\tau$ 
($\tau=\tau_3$) to be diagonal and label the states with the respective 
eigenvalues.

\beq
\begin{array}{lcl}
\vec{T}^2 | j,m,r,t \ra &=& [j][j+1] | j,m,r,t \ra \\[3mm]
\tau | j,m,r,t \ra &=& q^{-4m} | j,m,r,t \ra \\[3mm]
X^0 | j,m,r,t \ra &=& t | j,m,r,t \ra \\[3mm]
X \ci X  | j,m,r,t \ra &=& r^2 | j,m,r,t \ra
\end{array}
\eeq
The well known representations of the $T$ algebra are given in the appendix.
As in the undeformed case, the $T$, $X^A$ algebra allows us to express the 
``vector''-operator $X^A$ through reduced matrix elements \cite{pauli}.

The respective algebra as it follows from chapter 1 is:

\bea
\tau X^3&=& X^3 \tau \no
\tau X^+&=& q^{-4} X^+ \tau \label{TX1} \\
\tau X^-&=& q^4 X^- \tau \nonumber
\eea
\bea
T^- X^3 &=& X^3  T^-  + q \sqrt{1 + q^2} X^-\no
T^+ X^- &=& q^2  X^- T^+  + q^{-1} \sqrt{1 + q^2} X^3 \label{TX3} \\
T^- X^- &=& q^2 X^- T^-\nonumber
\eea
\bea
T^+ X^3 &=& X^3 T^+  + q^{-2} \sqrt{1 + q^2} X^+\no
T^+ X^+ &=& q^{-2} X^+ T^+\label{TX2} \\
T^- X^+ &=& q^{-2}  X^+ T^-  + \sqrt{1 + q^2} X^3 \nonumber
\eea

We proceed exactly as in the undeformed case.
{From} (\ref{TX1}) follows that $X^3$ does not change the eigenvalue of $\tau$
and that $X^+$, ($X^-$) changes $m$ by $+1$, ($-1$).

{From} (\ref{TX2}) we learn that the $m$ dependence of the $X^+$ matrix 
elements can be computed explicitly and that the matrix elements of $X^3$ 
can be expressed in terms of the reduced matrix elements of $X^+$.
{From} (\ref{TX3}) follow the same relations for $X^-$. Via $X^3$, the
reduced matrix elements of $X^+$ are related to the reduced matrix elements of 
$X^-$.
As $X^A$ commutes with $X \ci X$ and $X^0$, $X^A$ does not change the 
eigenvalues of $X^0$ and $X \ci X$.

For the non-vanishing matrix elements we obtain the following result:
\begin{eqnarray*}
\la j, m+1,r,t | X^+ | j, m, r,t  \ra &=&-q^{m+2} \sqrt{[j+m+1][j-m]} 
\la j,r,t  \| X^- \| j, r,t \ra \\
\la j+1, m+1,r,t | X^+ | j, m, r,t\ra &=&q^{m-2j} 
\sqrt{[j+m+1][j+m+2]} \\
&&\cdot \la j+1, r,t  \| X^- \| j, r,t  \ra  \\
\la j-1, m+1,r,t | X^+ | j, m, r,t \ra &=&q^{m+2j+2}
\sqrt{[j-m][j-m-1]} \\
&& \cdot \la j-1,r,t \| X^- \| j, r,t \ra
\end{eqnarray*}
\bea
\la j, m-1, r, t | X^- | j, m, r, t\ra &=&q^m
\sqrt{[j+m][j-m+1]} \la j, r, t \| X^- \| j, r, t \ra  \no
\la j+1, m-1, r, t | X^- | j, m, r, t\ra &=&q^m
\sqrt{[j-m+1][j-m+2]} \label{redXm}\\
&&\cdot \la j+1, r, t \| X^- \| j, r, t \ra  \no
\la j-1, m-1, r, t | X^- | j, m, r, t\ra &=&q^m
\sqrt{[j+m][j+m-1]} \la j-1, r, t \| X^- \| j, r, t \ra  \nonumber
\eea
\begin{eqnarray*}
\la j, m, r, t | X^3 | j, m, r, t\ra &=&q^{3 \over 2} 
{{\sqrt{1+q^2}}\over{q^2-1}} \{q^{2m}- {{q^{2j+1}+q^{-2j-1}}
\over{q+q^{-1}}}\} \la j, r, t \| X^- \| j, r, t \ra  \\
\la j+1, m, r, t | X^3 | j, m, r, t\ra &=&q^{m-j-{1 \over 2}}
\sqrt{1+q^2} \sqrt{[j-m+1][j+m+1]} \\
&& \cdot \la j+1, r, t \| X^- \| j, r, t \ra \no
\la j-1, m, r, t | X^3 | j, m, r, t\ra &=&-q^{m+j+{1 \over 2}}
\sqrt{1+q^2} \sqrt{[j-m][j+m]} \\
&& \cdot \la j-1, r, t \| X^- \| j, r, t \ra
\end{eqnarray*}

All the $m$ dependence of the $X^A$ matrix elements is now explicitly
known.

To get information on the reduced matrix elements we have to use the $X,X$
relations (\ref{XX}).

We start with the relation:
\beq
X^3 X^+-q^2X^+ X^3=(1-q^2) X^0 X^+
\label{X3Xp}
\eeq
Depending on what matrix elements we take, (\ref{X3Xp}) leads to a recursion 
formula for \linebreak
$\la j, r, t \| X^- \| j, r, t \ra$ or for the quantity $\rho_{r,t}$ 
which is defined as follows:
\beq
\rho_{r,t}(j+1)=[2j+1][2j+3] \la j,r,t\|X^-\|j+1,r,t \ra 
\la j+1,r,t\|X^-\|j,r,t\ra 
\label{rho}
\eeq
These recursion formulas can be solved and we obtain:
\beq
\la j, r, t \| X^- \| j, r, t \ra=-\lambda q^{-1} \sqrt{[2]}
{[j][j+1] \over [2j][2j+2]} t
\label{jXj}
\eeq
As $X^-$ changes the eigenvalue of $\tau$ the above matrix element
has to be zero for $j=0$. Eqn. (\ref{jXj}) is valid for $j \ge 1$.
\beq
\rho_{r,t}(j+1)=\rho_{r,t}(1)+\lambda^2 [2] q^{-2} t^2 
\sum_{l=1}^j {[2l+1] \over \{l\}^2 \{l+1\}^2}
\label{rho2}
\eeq
The quantity $\rho_{r,t}(1)$ is the unknown left. 
It is related to the radius $r$.
To see this we decompose $X \ci X$ into the product of matrix elements
of $X^A$. The calculation is particularly simple for the matrix element:
\beq
\la 0,0,r,t|X \ci X|0,0,r,t \ra=r^2=-q^2 [2] \rho_{r,t}(1)
\label{r1}
\eeq
For $j \neq 0$ the same calculation, but more tedious, yields:
\beq
\rho_{r,t}(j+1) ={1 \over q^2 [2]} \left\{-r^2+{[j][j+2]\lambda^2 \over 
\{j+1\}^2} t^2 \right\}
\label{rhofin}
\eeq
The two formulas (\ref{rho2}) and (\ref{rhofin}) agree because the sum
in (\ref{rho2}) can be summed up:
\beq
\sum_{l=1}^j {[2l+1] \over \{l\})^2\{l+1\}^2}=
{[2j] \over
[2] \{j\}^2\{j+1\}^2}(1+{[2j+2] \over [2]})
\label{sum}
\eeq
Eqn. (\ref{sum}) can be proved by induction.

From (\ref{r1}) follows that $\rho_{r,t}(1)$ is negative. We shall show that
this is true for $\rho_{r,t}(j)$ for any value of $j$. 
It follows from the hermiticity properties of the coordinates:
\beq
\overline{X^3}=X^3, \qquad  \overline{X^+}=-qX^-,\qquad  
\overline{X^-}=-\iq X^+  
\eeq
For the reduced matrix elements this implies:
\bea
\overline{\la j,r,t \| X^- \| j,r',t'\ra} &=& 
\la j,r',t' \| X^- \| j,r,t\ra \no
\overline{\la j+1,r,t \| X^- \| j,r',t'\ra} &=& -q^{2(j+1)}
\la j,r',t' \| X^- \| j+1,r,t\ra \label{Xmatcon}
\eea
{From} the definition of $\rho_{r,t}$ (\ref{rho}) we now obtain:
\beq
\rho_{r,t}(j+1)=-q^{-2(j+1)}[2j+1][2j+3] |\la j+1,r,t\|X^-\|j,r,t\ra |^2
\label{rho3}
\eeq
Thus $\rho_{r,t}(j)$ is negative or zero. 
This can lead to an upper bound for $j$.

If we combine (\ref{rho3}) with (\ref{rhofin}) we see that 
$|\la j+1,r,t\|X^-\|j+1,r,t \ra |^2$
is now explicitly known as a function of $j$, $r$ and $t$. Only a phase is 
undetermined.
But the relative phase between states with different $j$ has not been fixed 
yet.
We do it now, by assuming that $\la j+1,r,t\|X^-\|j,r,t \ra$
is real. Eqn. (\ref{Xmatcon}) then determines $\la j,r,t\|X^-\|j+1,r,t \ra$.

As the reduced matrix elements $\la j,r,t\|X^-\|j,r,t \ra$ have been
given in (\ref{jXj}), all the matrix elements of $X^A$ are known
as functions of $j$, $m$, $r$ and $t$.

\vspace{1ex}
We have to learn more about the spectrum of $t$ and $r$. This can be done by
studying the $X$, $R$ algebra.

We start with the following relations, they are a consequence of our algebra: 
\bea
U X^0&=&{1 \over q} {{q^4+1} \over {q^2+1}} X^0 U-{1 \over q} (q^2-1)^2
X \ci R \label{algebra} \\
X \ci R \; X^0&=&{2q \over 1+q^ 2} X^0 \; X \ci R -{q \over (1+q^2)^2} 
X \ci X \; U
\nonumber
\eea
If we take matrix elements of these relations we get two homogeneous linear 
equations in the matrix elements of $U$ and $ X \ci R$ that have a non trivial 
solution only if the determinant of the coefficient matrix vanishes:
\beq
\left(t-{\lambda \over [2]} t' \right) \left(t-{\{2\} \over [2]} t' \right) 
-{\lambda^2 \over [2]^2} {r'}^2=0
\label{det}
\eeq
The invariant length commutes with $X \ci R$ and $U$, and, as a consequence
\beq
s^2=t^2-r^2={s'}^2={t'}^2-{r'}^2=-l^2
\label{sinv}
\eeq
We shall use as a variable $s^2$ for the time-like and $l^2=-s^2$
for the space-like case.

If we replace ${r'}^2$ in (\ref{det}) by $s^2$ we obtain a quadratic 
equation in $t$ that has the solution:
\beq
t={[2] \over 2}t'\pm {\lambda \over 2} \sqrt{{t'}^2
-{\left(2 \over [2] \right)^2} s^2}
\label{eigt}
\eeq
Thus $t$ and $t'$ have to be related this way for a non-vanishing matrix
element $X \ci R$.
For ${r'}^2=0$ however there is a special situation. 
From (\ref{rhofin}) follows that $\rho_{0,t}(1)=0$ and that $j$ has to be 
zero for  $\rho_{0,t}$ not to be positive.
From (\ref{rho3}) all the $X$ matrix elements are zero and thus the 
$X \ci R$ matrix elements as well:
\beq
\la 0,0,0,t' | X \ci R | 0,0,r,t \ra=0
\eeq
The second equation of (\ref{algebra}) is trivially satisfied leaving us 
with the equation:
\beq
\la 0,0,0,t' | U | 0,0,r,t \ra \left(t-{\{2\} \over [2]} t'\right)=0
\eeq
In this case the matrix element of $U$ can only be different from zero if
\beq
t={\{2\} \over [2]} t'
\eeq
In all the other cases, the matrix element of $U$ is related to the matrix 
element of $X \ci R$:
\beq
\la j,m,r',t' |U | j,m,r,t \ra \left(t- {\{2\} \over [2]} t' \right)=
-q \lambda^2 \la j,m,r',t' |X \ci R | j,m,r,t \ra
\label{U1}
\eeq

We shall now discuss the time-like, space-like and light-like region
separately.

Let us start with the time-like region $s^2>0$. We assume that there is a
point $r'=0$ on the hyperbola, thus $t'_0=s$ ($t'_0=-s$).

According to the discussion above there is only the matrix element of $U$ that
connects this state to the state to the time:
\beq
t_1={\{2\} \over [2]} s > s
\eeq
We now continue to use  (\ref{eigt}) and find the other values of $t$:
\beq
t_n={s \over [2]} \{n+1\}  \qquad n=0,1,\ldots,\infty
\label{ttime}
\eeq
The values for $r_n$ follow from (\ref{sinv}):
\beq
r_n^2={s^2 \lambda^2 \over [2]^2} [n+2][n]
\label{rtime}
\eeq
For the backward light-cone we just have to take $t_0$ negative.

If we would not have assumed $r'=0$ to be in the spectrum our matrix elements
would connect to negative values of $r^2$.

For the space-like region, $s^2=-l^2<0$, a similar analysis gives the following
values for $t$ and $r$:
\bea
t_n&=&\pm {l \lambda \over [2]}[n] \qquad n=-\infty \ldots \infty 
\label{trspace} \\
r_n^2&=&{l^2 \over [2]^2} \{n+1\}\{n-1\}
\nonumber
\eea
On the light cone, $s^2=0$:
\bea
t_n&=&q^n \tau_0 \qquad n=-\infty \ldots \infty \\
r^2_n&=&q^{2n} \tau_0^2 \nonumber
\eea

If we now go back to eqn. (\ref{rhofin}) and insert the values of $r_n$ and
$t_n$ we find that for the time-like region $\rho(n+1)=0$.
That means that in this case $j$ is restricted to be $j \le n$.
There is no restriction of this type for the space-like region of the 
light-cone.

To conclude this section we give an explicit form for $\rho$ for the 
time-like, space-like and light-like region. 
We find:
\beq
\begin{array}{llcl}
\mbox{space-like:} & \quad \rho_n(j+1)&=&-{\displaystyle l^2 \over 
\displaystyle [2] q^2} {\displaystyle \{n-j-1\}\{n+j+1\} \over 
\displaystyle \{j+1\}^2} \hspace*{4cm}\\[4mm]
\mbox{time-like:} & \quad \rho_n(j+1)&=&-{\displaystyle s^2 \lambda \over
\displaystyle [2] q^2} {\displaystyle [n-j-1][n+j+2] \over \displaystyle 
\{j+1\}^2} \\[4mm]
\mbox{light-like:} & \quad \rho_n(j+1)&=&-{\displaystyle [2] \tau_0^2 
\over \displaystyle q^2} {\displaystyle 1 \over \displaystyle \{j+1\}^2}
\end{array}
\eeq
We see that only for the time-like region $\rho_n$ can change sign.

{From} (\ref{rho3}) follows with our phase convention:
\bea 
\la j+1,r,t \| X^- \| j,r,t\ra &=&q^{j+1} \sqrt{-\rho_{r,t}(j+1) \over
[2j+1][2j+3]} \\
\la j,r,t \| X^- \| j+1,r,t\ra &=&-q^{-j-1} \sqrt{-\rho_{r,t}(j+1) \over
[2j+1][2j+3]} \nonumber
\eea
We already know $\la j,r,t \| X^- \| j,r,t \ra$ (\ref{jXj}).
Then all the matrix elements of $X$ depend on $s$ for the time-like,
on $l$ for the space-like and on $\tau_0$ for the light-like
region as the only undetermined variable.

\Section{Matrix elements of the generators \boldmath $R^A$ of the \\
\boldmath $q$-Lorentz algebra}

The operators $R^A$ are ``vector'' operators as well, and their matrix elements
can be expressed through reduced matrix elements. The formulas (\ref{redXm})
are valid for $R^A$ except that $R^A$ is not diagonal in $r$
and $t$. 

If we analyze the ``scalar'' product of two arbitrary ``vector'' operators
through  matrix elements we get the general formula:
\bea
\la j,m,\mu | A \ci B | j,m,\nu\ra &=&\sum_{\nu'} 
\la j,\mu \| A^- \| j,\nu'\ra 
\la j,\nu' \| B^- \| j,\nu\ra
{1 \over [2]} q^2 [2j+2][2j]\label{matprod} \\
&&-\la j,\mu \| A^- \| j+1,\nu'\ra \la j+1,\nu' \| B^- 
\| j,\nu \ra
q^2 [2j+2][2j+3] \no
&&- \la j,\mu \| A^- \| j-1,\nu'\ra \la j-1,\nu' \| B^- 
\| j,\nu \ra
q^2 [2j] [2j-1] \nonumber 
\eea
where $\nu,\mu$ stand for the quantum numbers $t$ and $r$.

We can apply (\ref{matprod}) to $X \ci R$ and find:
\bea
\la j,m,\mu | X \ci R | j,m,\nu\ra &=& 
\la j,\mu \| X^- \| j,\mu \ra 
\la j,\mu \| R^- \| j,\nu\ra
{1 \over [2]} q^2 [2j+2][2j] \label{XciRfin} \\
&&-\la j,\mu \| X^- \| j+1,\mu \ra \la j+1,\mu \| R^- 
\| j,\nu \ra
q^2 [2j+2][2j+3] \no
&&- \la j,\mu \| X^- \| j-1,\mu \ra \la j-1,\mu \| R^- 
\| j,\nu \ra
q^2 [2j] [2j-1] \nonumber 
\eea
This shows that in general ($r_\mu \neq 0$) non vanishing matrix elements 
of $R^A$ will lead to non vanishing matrix elements of $X \ci R$. We know from
the last section that $X \ci R$ has only non vanishing matrix elements 
between states labeled by $t_n$ and $t_{n \pm 1}$. Thus the 
non vanishing matrix elements for $R^A$ are between these states as well.

We now use the $R$,$X$ relations to get information on the $R$ matrix elements.
First the algebraic relation:
\beq
R^+ X^+=q X^+ R^+
\label{RpXp}
\eeq
If we take the $(j+2,j)$ matrix elements of this equation we obtain the 
recursion formulas:
\beq
{{\la j+2, r', t' \| R^- \| j+1,r, t\ra } \over
{\la j+1, r',t' \| R^- \| j,r, t\ra }}=
q {{\la j+2,r', t' \| X^- \| j+1,r',t'\ra } \over
{\la j+1, r,t \| X^- \| j, r,t\ra }}
\label{recR1}
\eeq
and:
\beq
{{\la j-2, r', t', \| R^- \| j-1, r, t\ra } \over 
{\la j-1, r',t' \| R^- \| j, r, t\ra }}=
q {{\la j-2, r', t' \| X^- \| j-1, r', t'\ra } \over 
{\la j-1, r, t \| X^- \| j, r, t\ra }}
\label{recR2}
\eeq
These formulas can be iterated. With the matrix elements of $X^-$ expressed 
in terms of $\rho$, we find for $j > 1$:
\beq
{\la j+1, r',t' \| R^- \| j,r, t\ra \over \la 1, r',t' \| R^- \| 0,r, t\ra}
=q^{2j} \sqrt{[3] \over [2j+1][2j+3]} \sqrt{\rho_{r',t'}(j+1) 
\ldots \rho_{r',t'}(2) \over \rho_{r,t}(j) \ldots \rho_{r,t}(1)}
\label{jpRj}
\eeq
and
\[
{\la j, r',t' \| R^- \| j+1,r, t\ra \over \la 0, r',t' \| R^- \| 1,r, t\ra}
=q^{-2j} \sqrt{[3] \over [2j+1][2j+3]} \sqrt{\rho_{r,t}(j+1) 
\ldots \rho_{r,t}(2) \over \rho_{r',t'}(j) \ldots \rho_{r',t'}(1)}
\]
There is another relation that follows from (\ref{RpXp}) 
if we take the $(j+1,j)$ matrix elements. It is:
\bea
&\la j+1, r', t' \| R^- \| j+1,r, t\ra 
\la j+1, r,t \| X^- \| j, r,t\ra \no
&+\la j+1, r',t' \| R^- \| j,r, t\ra
\la j, r,t \| X^- \| j, r,t\ra= \label{RXm2}\\
&q \la j+1,r', t' \| X^- \| j+1, r',t'\ra
\la j+1, r',t' \| R^- \| j,r, t\ra \no
&+q \la j+1,r', t' \| X^- \| j, r',t'\ra
\la j, r',t' \| R^- \| j,r, t\ra \nonumber
\eea
This equation is valid for $j \ge 1$ and relates ($j+1,j+1$),($j,j$) and 
($j+1,j$) matrix elements of $R^-$.

If we study the relation:
\beq
R^+(X^3-X^0)=\iq (X^3-X^0) R^+
\eeq
and its ($j+1,j$) matrix elements, the same $R^-$ matrix elements
as (\ref{RXm2}) are related. They can be combined to eliminate the ($j+1,j+1$)
matrix elements and to give a relation between 
the ($j,j$) and ($j+1,j$) matrix elements of $R^-$:
\bea
&\la j, r',t' \| R^- \| j,r, t\ra \la j+1, r',t' \| X^- \| j,r', t'\ra q^{j+2}
\\
&=\la j+1, r',t' \| R^- \| j,r, t\ra \sqrt{[2]} 
\left({\displaystyle t'\{j\}-t\{j+1\}} \over {\displaystyle \{j\}\{j+1\}}
\right)
\nonumber
\eea
It is valid for $j \ge 1$.

Taking the corresponding ($j-1,j$) matrix elements we obtain:
\bea
&\la j, r',t' \| R^- \| j,r, t\ra \la j-1, r',t' \| X^- \| j,r', t'\ra q^{-j+1}
\\
&=\la j-1, r',t' \| R^- \| j,r, t\ra \sqrt{[2]} 
\left({\displaystyle t'\{j+1\}-t\{j\}}
\over {\displaystyle \{j\}\{j+1\}} \right)
\nonumber
\eea
This equation is valid for $j > 1$.

Both equations can be used to find the ($j,j$) matrix elements
from (\ref{jpRj}) in terms of the $(1,0)$ or $(0,1)$ matrix elements of $R^-$.
\bea
\la j, r',t' \| R^- \| j,r, t\ra&=&\la 1, r',t' \| R^- \| 0,r, t\ra {1 \over r}
{t'\{j\}-t\{j+1\} \over \{j\}\{j+1\}}{[2] \sqrt{[3]} \over  q^2} \no
&& \cdot \sqrt{\rho_{r',t'} (j) \ldots \rho_{r',t'}(2) \over \rho_{r,t} (j) 
\ldots \rho_{r,t}(2)} \label{jRj}\\
&=&-\la 0, r',t' \| R^- \| 1,r, t\ra {1 \over r'}
{t'\{j+1\}-t\{j\} \over \{j\}\{j+1\}}q^2 [2] \sqrt{[3]} \no
&&\cdot \sqrt{\rho_{r,t} (j) \ldots \rho_{r,t}(2) \over \rho_{r',t'} (j) 
\ldots \rho_{r',t'}(2)}
\nonumber
\eea
For the values of $\rho$ given in (\ref{rhofin}) it can be seen by induction
in $j$ that the relation between $\la 1, r', t' \| R^- \| 0,r, t\ra$ and 
$\la 0, r', t' \| R^- \| 1,r, t\ra$ that follows from (\ref{jRj}) is
indeed independent of $j$. We take $j=2$ and obtain:
\beq
\la 1, r',t' \| R^- \| 0,r, t\ra=\la 0, r',t' \| R^- \| 1,r, t\ra (-q^4)
{r \over r'} {t' \{3\} -t \{2\} \over t' \{2\} -t \{3\}} {\rho_{r,t} (2) \over
\rho_{r',t'}(2)} 
\eeq
One of the matrix elements, e.g. $\la 1, r', t' \| R^- \| 0,r, t\ra$,
remains to be determined.

We already know  that the $U$ matrix elements are related to the $R^-$
matrix elements  from (\ref{U1}) and that $U$ is hermitean.
\beq
\la 0,0,t,r | U | 0,0,t',r'\ra=\overline{\la 0,0,t',r' | U | 0,0,t,r\ra}
\eeq
If we now use the relation:
\beq
U R^+ =R^+U
\label{URp}
\eeq
we find:
\beq
\Gamma(n)=|\la 0,0,t_n,r_n | U | 0,0,t_{n+1},r_{n+1}\ra|^2=\Gamma(n+1)
\eeq
Thus $\Gamma(n)$ is $n$ independent. This is valid for the time-like, 
space-like and light-like regions. 

To finally determine $\Gamma(n)$ we have to use a relation that fixes 
the length of $R \ci R$. This relation is:
\beq
U^2-1=(q^4-1)^2 R \ci R
\label{cas}
\eeq
This is now sufficient to determine $\Gamma(n)$. We find:
\beq
\Gamma(n)={1 \over [2]^2}
\eeq

We can use the freedom of choosing the phase of states with different
$n$ eigenvalues such that:
\beq
\la 0,0,r,t | U | 0,0,r',t'\ra={1 \over [2]}
\eeq

This determines all the matrix elements of $U$, $R^A$ and $S^A$, as
the $S^A$ matrix elements are conjugate to the $R^A$ matrix elements
(\ref{star}).

We finally give the explicit form of the following $R^-$ matrix elements:
\pn For $s^2=0$
\beq
\begin{array}{lcl}
\la 1,r_n,t_n \| R^- \| 0,r_{n-1},t_{n-1} \ra&=& {1 \over [2]^{\frac{5}{2}} 
\sqrt{[3]} \lambda} \\[4mm]
\la 1,r_n,t_n \| R^- \| 0,r_{n+1},t_{n+1} \ra&=&-\iqt {1 \over [2]^{\frac{5}{2}}
\sqrt{[3]} \lambda} \\[4mm]
\la 0,r_n,t_n \|R^-\|1,r_{n-1},t_{n-1} \ra &=&-{1 \over  [2]^{\frac{5}{2}} 
\sqrt{[3]} q^6 \lambda}  \\[4mm]
\la 0,r_n,t_n \|R^-\|1,r_{n+1},t_{n+1} \ra &=&{1 \over  [2]^{\frac{5}{2}} 
\sqrt{[3]} q^4 \lambda}
\end{array}
\eeq  
\pn For $s^2$ time-like
\beq
\begin{array}{lcl}
\la 1,r_n,t_n \| R^- \| 0,r_{n-1},t_{n-1} \ra&=&
{1 \over [2]^{\frac{5}{2}} \sqrt{[3]} q \lambda} \sqrt{[n+2] \over [n]} \\[4mm]
\la 1,r_n,t_n \| R^- \| 0,r_{n+1},t_{n+1} \ra&=&-{1 \over [2]^{\frac{5}{2}} 
\sqrt{[3]} q \lambda} \sqrt{[n] \over [n+2]} \\[4mm]
\la 0,r_n,t_n \|R^-\|1,r_{n-1},t_{n-1} \ra &=&-{1 \over  [2]^{\frac{5}{2}} 
\sqrt{[3]} q^5 \lambda} \sqrt{[n-1] \over [n+1]} \\[4mm]
\la 0,r_n,t_n \|R^-\|1,r_{n+1},t_{n+1} \ra &=&{1 \over  [2]^{\frac{5}{2}} 
\sqrt{[3]} q^5 \lambda} \sqrt{[n+3] \over [n+1]}
\end{array}
\eeq  
\pn For $s^2$ space-like
\beq
\begin{array}{lcl}
\la 1,r_n,t_n \| R^- \| 0,r_{n-1},t_{n-1} \ra&=&{1 \over  [2]^{\frac{5}{2}} 
\sqrt{[3]} q \lambda} \sqrt{\{n+1\} \over \{n-1\}} \\[4mm]
\la 1,r_n,t_n \| R^- \| 0,r_{n+1},t_{n+1} \ra&=&-{1 \over [2]^{\frac{5}{2}} 
\sqrt{[3]} q \lambda} \sqrt{\{n-1\} \over \{n+1\}} \\[4mm]
\la 0,r_n,t_n\|R^-\|1,r_{n-1},t_{n-1} \ra &=&-{1 \over [2]^{\frac{5}{2}} 
\sqrt{[3]} q^5 \lambda} \sqrt{\{n-2\} \over \{n\}} \\[4mm]
\la 0,r_n,t_n \|R^-\|1,r_{n+1},t_{n+1} \ra &=&{1 \over  [2]^{\frac{5}{2}} 
\sqrt{[3]} q^5 \lambda} \sqrt{\{n+2\} \over \{n\}}
\end{array}
\eeq 

We could have started from the momenta instead of the coordinates,
then we would have constructed representations of the $q$-deformed
Poincar{\'e} algebra. Such representations are obtained by replacing
$X^a$ everywhere with $P^a$ \cite{poincare}.

It should be noted that the representations with positive mass square
${p^0}^2-\vec{p}^2 >0$, have angular momentum limited by $j \le n$
(see discussion after the eqn. (\ref{U1})).

\Section{The scaling operator \boldmath$\Lambda$ and the spectrum of 
\boldmath$X^0$, \boldmath$X \ci X$} 

The action of the scaling operator $\Lambda^{\frac{1}{2}}$ on the states 
$|j,m,r,t \ra$ is easily found from (\ref{la}):
\beq
\Lambda^{\frac{1}{2}} |j,m,r,t \ra = \alpha_{j,m,r,t} |j,m,qr,qt \ra    
\label{alfa}
\eeq
{From} (\ref{star}) and (\ref{alfa}) follows:
\beq
|\alpha_{j,m,r,t}|^2=q^4
\eeq
It is obvious that $\Lambda^{\frac{1}{2}}$ changes the value of $s^2$ by a 
factor $q^2$. This shows that the values of $s$ and $l$ in (\ref{ttime}),
(\ref{rtime}) and (\ref{trspace}) have to take the
following values:
\bea
&& s=t_0 q^M, \qquad M=-\infty \ldots \infty \\
&& l=l_0 q^M
\nonumber
\eea
It is only the light-cone that is left invariant under the action of 
$\Lambda^{\frac{1}{2}}$

The states can be labeled with $j$, $m$, $n$ and $M$ for $s^2 >0$ and for
$s^2<0$. For $s^2=0$, $j$, $m$ and $n$ are sufficient.

\pn For $s^2>0$:
\begin{eqnarray*}
M=-\infty \ldots \infty \\
n=0 \ldots \infty \\
j=0 \ldots n
\end{eqnarray*}
\bea
X^0 \; |j,m,n,M \ra &=& {t_0 q^M \over [2]} \{n+1\} \; |j,m,n,M \ra \\
X \ci X \; |j,m,n,M \ra &=&  {t_0^2 q^{2M} \lambda^2 \over [2]^2} [n+2] [n] 
\; |j,m,n,M \ra \no
\Lambda^{\frac{1}{2}}  \; |j,m,n,M \ra &=& q^2 |j,m,n,M+1 \ra \nonumber
\eea
\pn For $s^2<0$:
\begin{eqnarray*} 
M=-\infty \ldots \infty \\
n=-\infty \ldots \infty \\
j=0 \ldots \infty
\end{eqnarray*} 
\bea
X^0 \; |j,m,n,M \ra &=& {l_0 q^M \over [2]} \lambda [n] \; |j,m,n,M \ra \\
X \ci X \; |j,m,n,M \ra &=&  {l_0^2 q^{2M} \over [2]^2} \{n+1\} \{n-1\} \; 
|j,m,n,M \ra \no
\Lambda^{\frac{1}{2}}  \; |j,m,n,M \ra &=& q^2 |j,m,m,n,M+1 \ra \nonumber
\eea
\pn For $s^2=0$:
\begin{eqnarray*} 
n=-\infty \ldots \infty \\
j=0 \ldots \infty
\end{eqnarray*} 
\bea
X^0 \; |j,m,n \ra &=& \tau_0 q^n \; |j,m,n \ra \\
X \ci X \; |j,m,n \ra &=& \tau_0^2 q^{2n} \; |j,m,n \ra \no
\Lambda^{\frac{1}{2}}  \; |j,m,n \ra &=& e^{i \alpha_n} q^2 |j,m,m,n+1 \ra 
\nonumber
\eea
In this case we cannot use the freedom of phase for the states to have
$\alpha=0$. 

As we shall need the $U$-matrix elements in the next section we list them here
explicitly.

\pn For $s^2<0$:
\bea
&\la j,m,n,M | U |j,m,n+1,M \ra=\la j,m,n+1,M| U |j,m,n,M\ra=\\
&{\displaystyle {1 \over [2]} \sqrt{\{n-j\}\{n+j+1\} \over 
\{n\}\{n+1\}}} \nonumber
\eea

\pn For $s^2>0$:
\bea
&\la j,m,n,M | U |j,m,n+1,M \ra=\la j,m,n+1,M| U |j,m,n,M\ra=\\
&{\displaystyle {1 \over [2]} \sqrt{[n-j+1][n+j+2] \over 
[n+1][n+2]}} \nonumber
\eea
We see that for the time-like region the matrix element of $U$ is zero
for $n=j-1$. This is in agreement with the condition $n \ge j$.

\pn For $s^2=0$:
\beq
\la j,m,n | U |j,m,n+1 \ra=\la j,m,n+1,M| U |j,m,n,M\ra=
{\displaystyle {1 \over [2]}}
\label{Ulight}
\eeq

We shall see that for $s^2 \neq 0$ these states are sufficient to construct
a representation of the full algebra introduced in chapter 1.
For $s^2=0$ there is no representation of this algebra.

\Section{Matrix elements of the momenta}

We first write the $q$-deformed Heisenberg relations (\ref{PX}) in a more 
explicit version:
\bea
& q^2 [2] P^0 X^0-q \{2\} X^0 P^0 + \lambda X \ci  P \no
&={i \over 2}[2] \{2\} q^4 \Lambda^{-{1 \over 2}} U \label{XPe1a} \\
&q^2 [2] P^0 X^A- q \{2\} X^A P^0- \lambda q^2  X^0 P^A
-\lambda \epsilon_{DC}{}^A X^C P^D\no
&=-{i \over2} [2]^2 q^6 \lambda \Lambda^{-{1 \over 2}} (q^2 R^A+S^A)
\label{XPe1b} \\
&q^2 [2] P^A X^0- q \{2\} X^0 P^A - \lambda q^2 X^A P^0
-\lambda \epsilon_{DC}{}^A X^C P^D\no
&={i \over2} [2]^2 q^6 \lambda  \Lambda^{-{1 \over 2}} (R^A+q^2S^A) 
\label{XPe1c} \\
&[2] (P^AX^B-X^A P^B) +{2 \over q^3} \epsilon_{DC}{}^E \epsilon_E{}^{AB} 
X^C P^D \no
&+ {\lambda \over q^2} \left( g^{AB} X \ci  P - g^{AB}X^0P^0
+\epsilon_C{}^{AB} (X^C P^0+X^0 P^C) \right) \no
&=-{i \over2} [2] q^2 \Lambda^{-{1 \over 2}} \left( \{2\} 
g^{AB} U-q^2 \lambda [2] \epsilon_C{}^{AB} (R^C-S^C) \right) \label{XPe1d}
\eea

All these relations contain $P^0$, in the relations (\ref{XPe1c}) and 
(\ref{XPe1d}) $P^0$ is multiplied by $\lambda$.

We rearrange these relations to obtain two relations containing $P^0$
and $X \ci P$ as the only unknowns.

First we contract (\ref{XPe1d}) with $g_{AB}$:
\bea
&P \ci X- {1 \over q^3} {\{2\} \over [2] } X \ci P -q^4 \lambda {[3] \over [2]}
X^0 P^0 \label{XPe1e} \\
&=-{i \over2} q^2 [2] [3]  \Lambda^{-{1 \over 2}} U
\nonumber
\eea

Eqn. (\ref{XPe1a}) and (\ref{XPe1e}) together with their conjugates 
yield three independent equations:
\bea
P^0 X^0-X^0 P^0&=
&{i \over 2} (q^4 \Lambda^{-{1 \over 2}} + \Lambda^{1 \over 2}) U 
\label{PX2a} \\
P \ci X -X \ci  P&=
&-{i \over 2} [3] (q^4 \Lambda^{-{1 \over 2}} + \Lambda^{1 \over 2}) 
U \label{PX2b} \\
\lambda (X \ci P - X^0P^0)&=&
{i \over 2} q^2 [2] (\Lambda^{1 \over 2} -\Lambda^{-{1 \over 2}})U \label{PX2c}
\eea

Eqn. (\ref{PX2b}) can be used to express $P \ci X$ in terms of $X \ci P$.
Eqn. (\ref{PX2c}) is one of the wanted equations, the second one is obtained
by multiplying (\ref{XPe1c}) by $X^B g_{BA}$:

\beq
X \ci P X^0-{2 \over q [2]} X^0 X \ci P-{\lambda \over [2]} X \ci X P^0
=i q^4 \lambda [2] X \ci  R \Lambda^{-{1 \over 2}} \nonumber
\label{XPX0}
\eeq

This provides us with a system of two linear equations for the two unknowns,
the matrix elements of $P^0$ and $X \ci P$. The determinant of this
system of linear equations is proportional to $[2] t' (t-t')+ \lambda s^2$.
For $s^2 \neq 0$ the equations can be solved.
For $s^2 =0$ and $t=t'$ the determinant vanishes. The homogeneous
part of the two equations becomes linear dependent. For the inhomogeneous
part this would imply $\la j,m,n | U | j,m,n+1 \ra=0$,
in clear contradiction to (\ref{Ulight}).
We conclude that the $s^2=0$ representation of the $X^a$, $R^A$, $S^A$, $U$,
$\Lambda$ algebra cannot be extended to a representation of the full
algebra. For $s^2 \neq 0$ we can calculate the matrix elements. They are
consistent with (\ref{PX2a}) and the other algebra relations.

From the $X \ci P$ matrix elements we obtain the reduced matrix elements of
$P^-$, hermiticity of $P$ has to be used. This way we obtain all the matrix
elements of $P^A$.
Representations of the full algebra have now been constructed. Their
explicit form is given in chapter 2. It is interesting that the forward,
backward time-like and the space-like regions provide inequivalent, irreducible
representations by themselves.

%\newpage

\setcounter{section}{0}
\appendix{\boldmath $R$-matrices, metric and \boldmath $\varepsilon$-tensor}

\pn
{\bf Euclidean space}

\medskip

For the Euclidean space the metric tensor is defined as:
\begin{equation}
g_{AB} : \quad g_{+-} = - q,\quad  g_{33} = 1,\quad  g_{- +} = - \frac{1}{q}
\end{equation}
\[
g^{AB} : \quad g^{+-} = - q,\quad  g^{33} = 1,\quad  g^{- +} = - \frac{1}{q}
\]
\[
g_{AB} g^{BC}=\delta_A^C=g^{CB} g_{BA}
\]
With the metric indices can be raised and lowered:
\[
X_A = g_{AB} X^B, X^A = g^{AB} X_B
\]
and an invariant scalar product can be given:
\beq
X \circ Y = g_{AB} X^A Y^B=X^3Y^3-qX^+Y^--\frac{1}{q}X^-Y^+
\label{product}
\eeq
The $\varepsilon$-tensor is defined as:
\begin{eqnarray}
\varepsilon_{+-} {}^3 &=& q,\quad \varepsilon_{-+} {}^3 = -q,\quad
\varepsilon_{33} {}^3 = 1-q^2,\quad\no
\varepsilon_{+3} {}^+ &=& 1,\quad \varepsilon_{3+} {}^+ = -q^2,
\quad\nonumber\\
\varepsilon_{-3} {}^- &=& -q^2,\quad \varepsilon_{3-} {}^- =1.\label{epsilon}
\eea
Indices of the $\varepsilon$-tensor can also be raised and lowered through
the metric, e.g.:
\[
\varepsilon_{ABC} = g_{CD} \varepsilon_{AB}{}^D
\]
In terms of the metric and of the $\eps$-tensor the three-dimensional 
$\hat{R}$-matrix of the $q$-Euclidean space can be written in the form:
\beq
\hat{R}^{AB}{}_{CD} = \delta^A_C\delta^B_D -
q^{-4}\varepsilon^{FAB}\varepsilon_{FDC}-q^{-4} (q^2-1)
g^{AB} g_{CD}
\label{Rmatrexp}
\eeq

\medskip 

\pn
{\bf Minkowski space}

\medskip

For the $q$-deformed Minkowski space it turns out that two different 
$R$-matrices exist. Their projector decomposition is given by:
\begin{eqnarray}
\hat{R}_I & = & P_S+P_T-q^2P_+-q^{-2}P_- \label{R1} \\
& = & 1\hspace{-1mm}\mbox{\rm I}-(1+q^2)P_+-(1+\frac{1}{q^2})P_- \no
\hat{R}_{II} & = & q^{-2}P_S+q^2P_T-P_+-P_- \label{R2} \\
& = & \frac{1}{q^2}1\hspace{-1mm}\mbox{\rm I}+(q^2-\frac{1}{q^2})
P_T-(1+\frac{1}{q^2})P_A \nonumber
\end{eqnarray}
In these definitions $P_S$,$P_T$,$P_+$,$P_-$ are the projectors on the
symmetric, trace, selfdual, antiselfdual eigenspaces respectively.
This decomposition shows clearly that $\hat{R}_I$ cannot distinguish
the symmetric while $\hat{R}_{II}$ cannot distinguish the antisymmetric
eigenspaces, because they have the same eigenvalue, so that both matrices
are necessary to distinguish all the spaces. 
The explicit expression of the projectors follows.

\vspace{0.5cm}
\noindent $P_+$:

\beq
\begin{array}{c|cccc}
& \;\;00 & C0 & 0D & CD \\ \\
\hline \\
00 \;\;& \;\; 0 & 0 & 0 & 0  \\ \\
A0 \;\; & \;\;0 & \frac{q^2}{(1+q^2)^2} \delta_C^A & 
-\frac{1}{(1+q^2)^2}\delta_D^A &
\frac{1}{(1+q^2)^2}\varepsilon_{DC}{}^A \\ \\
0B \;\;&\;\; 0 & -\frac{q^4}{(1+q^2)^2}\delta_C^B & 
\frac{q^2}{(1+q^2)^2}\delta_D^B &
-\frac{q^2}{(1+q^2)^2}\varepsilon_{DC}{}^B \\ \\
AB \;\;& \;\;0 & \frac{q^2g^{EB}g^{FA}\varepsilon_{FEC}}{(1+q^2)^2} \quad & 
-\frac{g^{EB}g^{FA}\varepsilon_{FED}}{(1+q^2)^2} \quad  & 
\frac{\varepsilon_{DC} {}^E g^{SB} g^{RA} \varepsilon_{RSE}}{(1+q^2)^2} 
\quad \\
\end{array}
\eeq

\noindent $P_-$:

\beq
\begin{array}{c|cccc}
&\;\; 00 & C0 & 0D & CD \\ \\
\hline \\
00\;\; &\;\; 0   & 0 & 0 & 0 \\ \\
A0 \;\;& \;\;0 & \frac{q^2}{(1+q^2)^2} \delta_C^A & 
-\frac{q^4}{(1+q^2)^2}\delta_D^A &
-\frac{q^2}{(1+q^2)^2}\varepsilon_{DC}{}^A \\ \\
0B \;\;&\;\; 0 & -\frac{1}{(1+q^2)^2}\delta_C^B & 
\frac{q^2}{(1+q^2)^2}\delta_D^B &
\frac{1}{(1+q^2)^2}\varepsilon_{DC}{}^B \\ \\
AB\;\; &\;\; 0 & -\frac{g^{EB}g^{FA}\varepsilon_{FEC}}{(1+q^2)^2} \quad& 
\frac{q^2g^{EB}g^{FA}\varepsilon_{FED}}{(1+q^2)^2} \quad & 
\frac{\varepsilon_{DC} {}^E g^{SB} g^{RA} 
\varepsilon_{RSE}}{(1+q^2)^2} \quad \\
\end{array}
\eeq

%\newpage

\noindent $P_T$:

\beq
\begin{array}{c|cccc}
& 00 \;\;& \;\;C0\quad & 0D & CD \\ \\
\hline \\
00 \;\;& \;\;\frac{q^2}{(1+q^2)^2}   & 0\quad  & \;0 & 
-\frac{q^2}{(1+q^2)^2} g_{CD}  \\ \\
A0 \;\;& \;\;0 & 0\quad & \;0 & 0 \\ \\
0B \;\;& \;\;0 & 0\quad  & \;0 & 0 \\ \\
AB \;\;&  \;\;-\frac{q^2}{(1+q^2)^2}g^{AB} & 0\quad  & \;0 & 
\frac{q^2}{(1+q^2)^2}g^{AB} g_{CD} \\
\end{array}
\eeq
\vspace{0.5cm}

\noindent It holds:
\beq
1\hspace{-1mm}\mbox{\rm I}=P_S+P_T+P_++P_- 
\eeq
Using $P_T$ it is possible to construct a 4-dimensional metric:
\beq
\begin{array}{ll}
\eta_{00}=-1, \quad & \eta_{33}=1\\
\eta_{+-}=-q, \quad & \eta_{-+}=-\frac{1}{q} \\
\eta^{ab}=\eta_{ab}
\end{array}
\eeq
which enables to raise and lower indices:
\beq
X_A = \eta_{AB} X^B, X^A = \eta^{AB} X_B
\eeq
and to define an invariant scalar product in 4 dimensions:
\bea
X \cdot Y& = & X^0Y^0-X^3Y^3+qX^+Y^-+\frac{1}{q}X^-Y^+ \\
& = & -\eta_{ab} X^a Y^b \nonumber
\eea
The sum of the selfdual and antiselfdual projectors defines the
$q$-deformed antisymmetrizer:
\beq
P_A=P_++P_-
\eeq
while their difference defines the $q$-deformed 4-dimensional 
$\varepsilon$-tensor:
\beq
\varepsilon^{ab}{}_{cd}=P_+^{ab}{}_{cd} -  P_-^{ab}{}_{cd}
\eeq

\end{document}